\documentclass[]{ifacconf}
\usepackage{amsfonts}

\usepackage{amsmath}
\usepackage{graphicx}
\usepackage[sort]{natbib}
\usepackage{color}
\usepackage{verbatim}

\usepackage{natbib}

\begin{document}
\begin{frontmatter}

\title{On Approximating Polynomial-Quadratic Regulator Problems}
\author[First]{Jeff Borggaard}
\author[First]{Lizette Zietsman}

\address[First]{Department of Mathematics, Virginia Tech, Blacksburg, VA, USA.}


\begin{abstract}
Feedback control problems involving autonomous polynomial systems are prevalent, yet there are limited algorithms and software for approximating their solution.  This paper represents a step forward by considering the special case of the regulator problem where the state equation has polynomial nonlinearity, control costs are quadratic, and the feedback control is approximated by low-degree polynomials.  As this represents the natural extension of the {\em linear-quadratic regulator} (LQR) and {\em quadratic-quadratic regulator} (QQR) problems, we denote this class as {\em polynomial-quadratic regulator} (PQR) problems.  
The present approach is amenable to feedback approximations with low degree polynomials and to problems of modest model dimension.  This setting can be achieved in many problems using modern model reduction methods.  The Al'Brekht algorithm, when applied to polynomial nonlinearities represented as Kronecker products leads to an elegant formulation.  The terms of the feedback control lead to large linear systems that can be effectively solved with an N-way generalization of the Bartels-Stewart algorithm.  We demonstrate our algorithm with numerical examples that include the Lorenz equations, a ring of van der Pol oscillators, and a discretized version of the Burgers equation.  The software described here is available on Github.
\end{abstract}


\end{frontmatter}

\section{Motivation}
Linear feedback control of autonomous nonlinear systems, such as those describing the behavior of fluids, can be sufficient to achieve stabilization--even for an unstable steady-state solution.  For example, this has been demonstrated through the rotational stabilization of the wake behind a circular cylinder, cf.~\cite{bergmann2005orc}, \cite{benner2016rsl}, \cite{borggaard2010lfc}, and \cite{borggaard2014mrd}.

There is a shortage of software tools for nonlinear problems in control and systems theory.  The general Matlab Nonlinear Systems Toolbox (NST) by \cite{krener2015NST} took a broad step toward delivering useful tools for a number of important problems.  Since we inherently encounter the {\em curse of dimensionality} in these problems, there is also a need to develop specialized tools for important classes of problems.  This paper addresses this by specifically solving the polynomial-quadratic regulator problem: minimizing a quadratic cost subject to a state equation with a polynomial nonlinearity. 

For example, linear feedback laws found by solving the linear-quadratic regulator (LQR) problem compute the linear feedback law as the solution to a single algebraic Riccati equation and have the property that the linear portion of the nonlinear system becomes stable~\citep{barbu2007les,barbu2001fip,barbu2006tbs,raymond2006fbs}.  Unfortunately, for nonlinear systems, this only guarantees local stability.  The ability of linear feedback to stabilize the steady-state solution depends on the initial condition, which must be sufficiently close to the steady-state.  An alternative would be to develop nonlinear feedback control laws that could offer the ability to  expand the radius of convergence (shown with a simple example in \cite{borggaard2018computation}).  However, these require us to approximate solutions to the Hamilton-Jacobi-Bellman (HJB) equations, e.g.~\cite{kunisch2004hjb,breiten2019feedback}.  The HJB equations are notoriously complex in the general case.  Nevertheless, if one considers the polynomial-quadratic regulator (PQR) problem, having autonomous state equations with polynomial nonlinearities and a quadratic control objective, there is sufficient structure in polynomial approximations to approximate solutions.  One strategy is to use state-dependent Riccati equations (SDRE), cf.~\cite{banks2007nonlinear,cimen2012survey,cloutier1997statedependent}, which requires a proper factorization of the problem~\citep{banks2007nonlinear}.  SDREs have recently been applied to incompressible flows in~\cite{benner2017nonlinear}.  Another strategy, that we pursue here, is to use polynomial approximations
to the HJB equation and associated feedback operator based on Al'Brekht's method \citep{navasca2000solution}.  As described in this paper, recasting Al'Brekht's method for autonomous polynomial systems using Kronecker products leads to computable polynomial feedback laws for modest problem sizes.\footnote{This was recently used for bilinear systems in \cite{breiten2017taylor}.}  The PQR problem also happens to be exactly what is needed to solve discretized versions of distributed parameter control problems where the nonlinearity is quadratic (such as the Navier-Stokes equations used as our motivation above).  This is particularly true when linear feedback laws are being based on LQR problems. As in the LQR case, suitable model reduction methods \citep{aubry1988dcs,holmes1996tcs,ahmad2015nim} are essential to forming a solution methodology for distributed parameter control problems with quadratic nonlinearities. First of all, the Riccati equation must be solved to compute the linear term, e.g. \cite{singler2008alr,benner2013esl,singler2016pod} and the curse-of-dimensionality still appears with higher-order polynomial approximations of the feedback law.  

In this paper, we briefly outline the HJB equations, the PQR problem, and polynomial approximations to the value function and the feedback control operators.  Our formulation leads to a sequence of linear systems in Kronecker product form after an initial solution to the algebraic Riccati equation.  A na\"ive construction of these matrices and other terms would quickly become prohibitive.  However, the structure lends itself to newly developed recursive tensor linear algebra that avoids assembly and other taxing of computer memory.  We present a numerical study with a set of control problems with quadratic and cubic state equations to investigate the advantages of using higher degree approximations of the optimal feedback control.

\section{Problem Formulation}
For any $t\geq 0$, let ${\bf x}(t)\in \mathbb{R}^n$
be the state variables, ${\bf u}(t)\in \mathbb{R}^m$ the control inputs, and ${\bf Q}\in\mathbb{R}^{n\times n}$ and ${\bf R}\in \mathbb{R}^{m\times m}$ be weighting matrices satisfying properties ${\bf Q}^T={\bf Q}\geq0$ and ${\bf R}^T={\bf R}> 0$.
The {\em running cost} is defined as the quadratic $\ell({\bf x},{\bf u}) \equiv {\bf x}^T{\bf Q}{\bf x}+{\bf u}^T{\bf R}{\bf u}$.
The polynomial-quadratic regulator problem then is to find a control ${\bf u}(\cdot)\in L_2(0,\infty;\mathbb{R}^{m})$ that solves
\begin{equation}
\label{eq:oc}
  \min_{\bf u} J({\bf x},{\bf u}) = \int_0^\infty \ell({\bf x}(t),{\bf u}(t))\ dt,
\end{equation}
subject to the system dynamics
\begin{equation}  \label{eq:full}
  \dot{{\bf x}}(t) = {\bf A}{\bf x}(t) + {\bf B} {\bf u}(t) + {\bf f}({\bf x}(t)),\quad {\bf x}(0) = {\bf x}_0,
\end{equation}
where 
${\bf A} \in\mathbb{R}^{n\times n} $ and ${\bf B} \in \mathbb{R}^{n\times m}$ are constant matrices and $ \mathbf{f}: \mathbb{R}^{n}\longrightarrow\mathbb{R}^{n}$ is a $p$-degree polynomial in the states.
%
%

We define the {\em value function} $v({\bf x}_0) = J({\bf x}^*(\, \cdot\, ;{\bf x}_0),{\bf u}^*(\cdot))$ to be the value of (\ref{eq:oc}) when the optimal control ${\bf u}^*$ and corresponding state ${\bf x}^*$ are found from the initial point ${\bf x}_0$.  The optimal control is given by the feedback relation
\begin{equation}
\label{eq:feedback}
  {\bf u}(t) = {\mathcal K}({\bf x}(t)),
\end{equation}
which satisfies the HJB partial differential equations
\begin{align}
\label{eq:HJB1}
  0 &= \frac{\partial v}{\partial {\bf x}}({\bf x})\left({\bf A}{\bf x} + {\bf B} {\cal K}({\bf x}) +  {\bf f}({\bf x}) \right)
       + \ell({\bf x},{\cal K}({\bf x})), \\
\label{eq:HJB2}
  0 &= \frac{\partial v}{\partial {\bf x}}({\bf x}) {\bf B}
       + \frac{\partial \ell}{\partial {\bf u}}({\bf x},{\cal K}({\bf x})).
\end{align}
Ideally, we could solve the HJB equations simultaneously for $v$ and ${\mathcal K}$, but this is not computationally feasible due
to the curse of dimensionality.  Therefore, it is natural to consider series solutions.  As suggested in \cite{krener2014series}
the algorithm proposed by Al'Brekht is effective and has mathematical justification.  We will show that this can be solved effectively when the polynomials are expressed in Kronecker product form.  We first review some useful features of Kronecker products.

\section{Notation and Properties}
Kronecker products have a rich history in the control literature, cf.~\cite{brewer1978kronecker} and \cite{simoncini2016computational}.
The Kronecker product of two matrices ${\bf X}\in \mathbb{R}^{i_x\times j_x}$ and ${\bf Y}\in \mathbb{R}^{i_y\times j_y}$, with entries $x_{ij}$ and $y_{ij}$, is defined as the block matrix ${\bf X}\otimes{\bf Y} \in \mathbb{R}^{i_xi_y \times j_xj_y}$ with entries
\begin{displaymath}
   {\bf X} \otimes {\bf Y} \equiv \left[ \begin{array}{cccc}
     x_{11}{\bf Y} & x_{12}{\bf Y} & \cdots & x_{1j_x}{\bf Y} \\
     x_{21}{\bf Y} & x_{22}{\bf Y} & \cdots & x_{2j_x}{\bf Y} \\
     \vdots  &         &        & \vdots    \\
     x_{i_x1}{\bf Y}& x_{i_x2}{\bf Y}& \cdots & x_{i_xj_x}{\bf Y} \end{array} \right].
\end{displaymath}
The {\em vec} operation on any matrix ${\bf Y}$,  produces the vector ${\bf y}$ of length $i_y j_y$ with the $(i+(j-1)i_y)$-st entry being $y_{ij}$ (for $i=1$:$i_y$ and $j=1$:$j_y$) and is written ${\bf y} = {\rm vec}({\bf Y})$.
We will repeatedly utilize the following properties of Kronecker products: $({\bf C}\otimes{\bf D})({\bf E}\otimes{\bf F}) = ({\bf CE})\otimes({\bf DF})$ and $({\bf C}\otimes{\bf D})^T = {\bf C}^T \otimes {\bf D}^T$.
We will also take advantage of the Kronecker-vec relationship
\begin{displaymath}
  {\bf K}={\bf X}{\bf V}{\bf Y}^T \quad \mbox{leads to} \quad {\rm vec}({\bf K}) = ({\bf Y}\otimes{\bf X}){\rm vec}({\bf V})
\end{displaymath}
and the derivative of $c({\bf x}) = {\bf c}_2^T ({\bf x}\otimes{\bf x})$ in direction ${\bf f}$ as
\begin{displaymath}
  \frac{\partial c}{\partial {\bf x}} {\bf f} = {\bf c}_2^T ({\bf f}\otimes{\bf x} + {\bf x}\otimes{\bf f}).
\end{displaymath}

\section{The Polynomial-Quadratic Regulator}
We present the Kronecker product description of the polynomial-quadratic regulator (PQR) problem.  We specifically
express the polynomial nonlinearity as
\begin{displaymath}
  {\bf f}({\bf x}) \equiv {\bf N}_2 ({\bf x}\otimes{\bf x}) + \cdots + {\bf N}_p ({\bf x}\otimes \cdots \otimes{\bf x}).
\end{displaymath}  
where ${\bf N}_k \in \mathbb{R}^{n\times n^k}$ for $k=2,\ldots,p$.  By defining ${\bf q}_2 \equiv {\rm vec}({\bf Q})$ and 
${\bf r}_2 \equiv {\rm vec}({\bf R})$, we can rewrite
\begin{displaymath}
  \ell({\bf x},{\bf u}) = {\bf q}_2^T\left( {\bf x}\otimes{\bf x} \right) + {\bf r}_2^T\left( {\bf u}\otimes{\bf u} \right).
\end{displaymath}
Following Al'Brekht's approach, we now expand the value function and feedback operator as polynomials, then generate
equations by substitution of all of these polynomial expressions into (\ref{eq:HJB1})-(\ref{eq:HJB2}) and matching equal degree terms up to a desired approximation order.  Thus, we define
\begin{displaymath}
  v({\bf x}) = \underbrace{{\bf v}_2^T \left( {\bf x}\otimes{\bf x} \right)}_{v^{[2]}({\bf x})} + \underbrace{{\bf v}_3^T \left( {\bf x}\otimes{\bf x}\otimes{\bf x} \right) }_{v^{[3]}({\bf x})} + \cdots
\end{displaymath}
and
\begin{displaymath}
  \mathcal{K}({\bf x}) = \underbrace{{\bf k}_1 {\bf x}}_{{\bf k}^{[1]}({\bf x})} + \underbrace{{\bf k}_2 \left( {\bf x}\otimes{\bf x} \right)}_{{\bf k}^{[2]}({\bf x})} + \underbrace{{\bf k}_3 \left( {\bf x}\otimes{\bf x}\otimes{\bf x} \right)}_{{\bf k}^{[3]}({\bf x})} +\cdots,
\end{displaymath}
%
%
%
%
%
where ${\bf v}_d\in \mathbb{R}^{n^d\times 1}$ and ${\bf k}_d \in \mathbb{R}^{m\times n^d}$ are to be determined. 


Substituting the expansions for the value function $v$ and the feedback operator $\mathcal{K}$ into (\ref{eq:HJB1}), then collecting 
$O({\bf x}^2)$ terms and factoring, we have
\begin{align}\nonumber
  {\bf v}_2^T \left( ({\bf A}+{\bf B}{\bf k}_1)\otimes {\bf I}_n)\! +\! {\bf I}_n\otimes({\bf A}+{\bf B}{\bf k}_1) \right)({\bf x}\otimes{\bf x}) \\
  + {\bf q}_2^T({\bf x}\otimes{\bf x}) + {\bf r}_2^T({\bf k}_1\otimes{\bf k}_1)({\bf x}\otimes{\bf x}) = 0.\label{eq:ARE1}
\end{align}
Similarly, gathering $O({\bf x})$ terms from (\ref{eq:HJB2}) and using our assumptions on ${\bf R}$ leads to
\begin{equation}\label{eq:k1}
  {\bf v}_2^T({\bf B}\otimes{\bf I}_n) + {\bf r}_2^T ({\bf I}_m\otimes {\bf k}_1) = 0.
\end{equation} 
As we would expect, this is the LQR solution for the linear problem (ignoring ${\bf f}$ in (\ref{eq:oc})-(\ref{eq:full})), where we
have ${\bf v}_2 = {\rm vec}({\bf V}_2)$ and ${\bf V}_2$ solves the algebraic Riccati equation (ARE)
\begin{displaymath}
  {\bf A}^T {\bf V}_2 + {\bf V}_2{\bf A} - {\bf V}_2 {\bf B}{\bf R}^{-1}{\bf B}^T {\bf V}_2 + {\bf Q} = {\bf 0}.
\end{displaymath}
With ${\bf V}_2$ in hand, we can set ${\bf k}_1 = -{\bf R}^{-1}{\bf B}^T{\bf V}_2$.

Gathering $O({\bf x}^2)$ and $O({\bf x})$ terms from (\ref{eq:HJB1})-(\ref{eq:HJB2}) produces ${\bf v}_2$ and ${\bf k}_1$.
The Al'Brekht algorithm repeats this for successively higher degree terms with the following simplification.  Note that
gathering $O({\bf x}^{d+1})$ terms in (\ref{eq:HJB1}) to obtain an equation for ${\bf v}_{d+1}$ will produce the term
\begin{displaymath}
  {\bf v}_2^T(({\bf B}{\bf k}_d)\otimes {\bf I}_n) + {\bf r}_2^T({\bf k}_d\otimes{\bf k}_1).
\end{displaymath}
Factoring out $({\bf k}_d\otimes {\bf I}_n)$ from the right and using (\ref{eq:k1}), we clearly show the known result that this term
always vanishes.  This effectively decouples the calculation of the terms ${\bf v}_{d+1}$ and ${\bf k}_d$ for all values of $d>1$.
We now describe each of these calculations separately below.

\subsection{Coefficients of ${\bf v}_{d+1}$}
The degree three terms in (\ref{eq:HJB1}) can then be written using the definition ${\bf A}_c = {\bf A}+{\bf Bk}_1$ as
\begin{equation}\label{eq:d2}
\begin{split}
  \left({\bf A}_c \otimes {\bf I}_n \otimes {\bf I}_n + {\bf I}_n \otimes {\bf A}_c \otimes {\bf I}_n  +
  {\bf I}_n \otimes {\bf I}_n \otimes {\bf A}_c\right)^T {\bf v}_3 \\
  = -\left({\bf N}_2\otimes {\bf I}_n + {\bf I}_n\otimes {\bf N}_2\right)^T{\bf v}_2.
\end{split}
\end{equation}
As noted above, this is independent of ${\bf k}_2$.  

To write the equations from matching higher degree terms in a more compact way, we define the {\em N-way Lyapunov matrix}
or a special {\em Kronecker sum} matrix, cf.~\cite{benzi2015decay},
\begin{equation}\label{eq:kroneckerSum}
  \mathcal{L}_d({\bf X}) \equiv \underbrace{{\bf X}\! \otimes {\bf I}_n\! \otimes\! \cdots \otimes {\bf I}_n}_{d\ {\rm terms}} + \underbrace{{\bf I}_n\! \otimes\! {\bf X}\! \otimes {\bf I}_n\! \otimes\! \cdots \otimes {\bf I}_n}_{d\ {\rm terms}} + \cdots.
\end{equation} 
Then the calculation of ${\bf v}_3$ in equation (\ref{eq:d2}) follows from solving an equation of the form
\begin{equation}
\label{eq:degree3}
  \mathcal{L}_3({\bf A}_c^T) {\bf v}_3 = -\mathcal{L}_2({\bf N}_2^T) {\bf v}_2.
\end{equation}
Once we have ${\bf v}_3$, we can readily compute ${\bf k}_2$ as shown in Section~\ref{sec:k} below.  The other terms in the series expansion of the value function lead to equations that have a similar form.
All of the left-hand-sides are generically the same $\mathcal{L}_{d+1}({\bf A}_c^T) {\bf v}_{d+1}$.  However, 
the right-hand-sides of the equations gather more terms due to the ${\bf r}_2$ term in (\ref{eq:HJB1})
and the interactions of the previously computed nonlinear feedback terms with previously computed terms of the value function (that are known and moved to the right-hand-side).  This process is clarified by explicitly collecting the next two sets of terms for $v({\bf x})$ below.  For $O({\bf x}^4)$, we have
\begin{equation}
\label{eq:degree4}
  \mathcal{L}_4({\bf A}_c^T) {\bf v}_4 = -\mathcal{L}_3(({\bf B}{\bf k}_2+{\bf N}_2)^T) {\bf v}_3 - ({\bf k}_2^T \otimes
  {\bf k}_2^T ) {\bf r}_2,
\end{equation}
which can be solved for ${\bf v}_4$ once ${\bf k}_2$ is computed using the
solution ${\bf v}_3$ from (\ref{eq:degree3}), 
and
\vspace{-.1in}
\begin{eqnarray}\nonumber
  \mathcal{L}_5({\bf A}_c^T) {\bf v}_5 &=& -\mathcal{L}_4(({\bf B}{\bf k}_2+{\bf N}_2)^T) {\bf v}_4 
  -\mathcal{L}_3(({\bf B}{\bf k}_3+{\bf N}_3)^T){\bf v}_3 \\
  &&- ( {\bf k}_2^T\otimes {\bf k}_3^T + {\bf k}_3^T\otimes {\bf k}_2^T ) {\bf r}_2.
\label{eq:degree5}
\end{eqnarray}
Again, once we compute ${\bf k}_3$ from ${\bf v}_4$, we have everything we need to compute ${\bf v}_5$.

In general, while calculation of the coefficients ${\bf v}_{d}$ is described by large linear systems ($\mathcal{L}_{d}({\bf A}_c^T)\in\mathbb{R}^{n^{d}\times n^{d}}$), there is a great deal of structure and sparsity that can be exploited.  This will be discussed in Section~\ref{sec:linear}.   It is also
immediately obvious that {\em without} the nonlinear terms ${\bf N}_p$ in our state equation, the right-hand-side in (\ref{eq:degree3}) would vanish leading to ${\bf v}_3={\bf 0}$. The remaining equations for ${\bf v}_{d+1}$ would have homogeneous right-hand-sides and thus ${\bf v}_{d+1}={\bf 0}$ for $d=2$ and higher.  This is consistent with the LQR theory.  We also see
that even if ${\bf f}$ is quadratic, all terms in the series for $v$ and $\mathcal{K}$ could be non-zero.  Therefore, we are only computing approximations to the nonlinear feedback laws for non-trivial ${\bf f}$.

\subsection{Coefficients of ${\bf k}_{d}$\label{sec:k}}
We now turn our attention to using (\ref{eq:HJB2}) to calculate ${\bf k}_d$ from ${\bf v}_{d+1}$.  This is
again straight-forward using the specialized Kronecker sum operator (\ref{eq:kroneckerSum}),
\begin{equation}
\label{eq:k_d}
  {\bf k}_d = -\frac{1}{2}{\bf R}^{-1} \left(\mathcal{L}_{d+1}({\bf B}^T) {\bf v}_{d+1}\right)^T.
\end{equation}

\subsection{Computing Right-Hand-Side Vectors}
The assembly and solution of linear systems with the form 
\begin{equation}
\label{eq:generic}
  \mathcal{L}_{d+1}({\bf A}_c){\bf v}_{d+1} = {\bf c}
\end{equation}
is  only feasible for small values of $d$ and $n$.  
The advantage of the Kronecker product structure is that we can perform operations with Kronecker product matrices {\em without} actually forming the large block matrix.  The main issue that we deal with in this section is calculating ${\bf c}$, the terms on the right-hand-sides of e.g. (\ref{eq:degree3})--(\ref{eq:degree5}) or (\ref{eq:k_d}). Solution of the system (\ref{eq:generic}) is described in the next section. 

To calculate ${\bf c}$ for (\ref{eq:degree3})--(\ref{eq:k_d}) involves two types of terms.  The first involves the multiplication of a Kronecker form with a vector ${\bf r}_2$.  Recall, e.g.~\cite{brewer1978kronecker}, that 
\begin{equation}\label{eq:kron*v}
  ({\bf X}\otimes{\bf Y}){\bf r}_{2} = {\rm vec}({\bf Y}^T{\bf R}{\bf X}),
\end{equation}
where ${\bf R}$ has the appropriate dimensions and ${\bf r}_{2} = {\bf vec}({\bf R})$.  Therefore, the terms involving ${\bf r}_2$ only require
matrix multiplications and no assembly of the Kronecker product is required.

The second type of term are products of the Kronecker sum with a ${\bf v}_{d+1}$: $\mathcal{L}_{d+1}({\bf X}){\bf v}_{d+1}$.
Using the definition of (\ref{eq:kroneckerSum}), we have to calculate $d+1$ different multiplications 
of the Kronecker products with ${\bf v}_{d+1}$.  This is simplified using the associativity of the Kronecker product and writing
\begin{displaymath}
  {\bf I}_{n^\ell} = \underbrace{{\bf I}_n\otimes \cdots\otimes {\bf I}_n}_{\ell\ {\rm terms}}.
\end{displaymath}
The multiplications can be reduced to three different cases
\begin{displaymath}
  ({\bf X} \otimes {\bf I}_{n^d}){\bf v}_{d+1}, \ \ ({\bf I}_{n^{d-\ell}}\otimes {\bf X}\otimes {\bf I}_{n^{\ell}}){\bf v}_{d+1},
  \ \ \mbox{and} \ \ ({\bf I}_{n^d}\otimes {\bf X}){\bf v}_{d+1}.
\end{displaymath}
Here the relation (\ref{eq:kron*v}) and the associative law for Kronecker products are useful.  The first and last terms
above can be handled by the appropriate reshaping of ${\bf v}_{d+1}$ and multiplying with ${\bf X}$ (the multiplication by ${\bf I}_{n^d}$ is trivial).  The associative law allows us to handle all of the intermediate terms recursively as
\begin{align*}
   ({\bf I}_{n^{d-\ell}}\otimes {\bf X}\otimes {\bf I}_{n^{\ell}}){\bf v}_{d+1} &= (({\bf I}_{n^{d-\ell}}\otimes {\bf X})\otimes {\bf I}_{n^{\ell}}){\bf v}_{d+1}\\
       &= ({\bf I}_{n^{d-\ell}}\otimes ({\bf X}\otimes {\bf I}_{n^{\ell}})){\bf v}_{d+1}.
\end{align*}
The grouping can be done to maximize the size of the free identity matrix.

\subsection{Linear System Solutions\label{sec:linear}}
The Kronecker structure leads to larger systems (\ref{eq:generic}), but are now ameneble to modern high performance algorithms~\citep{kolda2009tensor,chen2019recursive,simoncini2016computational}.  Many of these algorithms,
e.g. \cite{chen2019recursive}, utilize a real Schur factorization of the matrix ${\bf A}_c$.
For this study, we used the {\em recursive} algorithms in \cite{chen2019recursive} for Laplace-like equations. Their software was 
trivially modified to take advantage of the fact that the same term ${\bf A}_c$ appears in every block and gave the
system exactly the form (\ref{eq:kroneckerSum}).

As an alternative, we have also developed a solver that generalizes the Bartels-Stewart algorithm to systems of the form
$\mathcal{L}_d({\bf A}_c){\bf v} = {\bf b}$.  As a preprocessing step, a Schur decomposition is performed on ${\bf A}_c$.
If ${\bf A}_c={\bf U}{\bf T}{\bf U}^*$, then we can apply this using the factorization property of the Kronecker product
\begin{displaymath}
  \mathcal{L}_d({\bf A}) = ({\bf U}\otimes\cdots\otimes{\bf U})\mathcal{L}_d({\bf T})({\bf U}\otimes\cdots\otimes{\bf U})^*.
\end{displaymath}
As with the Bartels-Stewart algorithm, we can work with upper triangular systems and also find the familiar solvability condition
in terms of not having eigenvalues of ${\bf T}$ (namely ${\bf A}_c$) reflected across the imaginary axis.  This cannot happen in our application 
as ${\bf A}_c$ is a stable matrix.  Note that instead of solving the upper triangular 
system directly, a block backward substitution algorithm allows us to take advantage of the sparsity pattern that arises with
$\mathcal{L}_d({\bf T})$ as well as fast matrix multiplications using the Kronecker-vec property.

\section{Numerical Results}\label{sec:numerical}
We present three sets of results.  The first is a controlled Lorenz system, the second is a ring of van der Pol oscillators,
and the final example is a discretized control problem involving
the one-dimensional Burgers equation.

\subsection{Controlled Lorenz Equations}
As a first test example, we consider the feedback control of the Lorenz equations where
\begin{displaymath}
  {\bf A} = \left[ \begin{array}{rrr} -10 & 10 & 0 \\ 28 & -1 & 0 \\ 0 & 0 & -8/3
  \end{array} \right],  \qquad {\bf B} = \left[ \begin{array}{r} 1 \\ 0 \\ 0 \end{array} \right],
\end{displaymath}
and the nonzero entries of ${\bf N}_2\in \mathbb{R}^{3\times 9}$ are
\begin{displaymath}
  {\bf N}_2(2,3)\!=\!{\bf N}_2(2,7)\!=\!-\frac{1}{2} \quad\mbox{and}\quad {\bf N}_2(3,2)\!=\!{\bf N}_2(3,4)\!=\!\frac{1}{2}, 
\end{displaymath}
accounting for the $-x_1x_3$ term in the second equation and the $+x_1x_2$ term in the third, respectively.
We choose ${\bf Q}={\bf I}_3$ and ${\bf R}={\bf I}_1$ as control weights.  The solution of the open- and closed-loop
systems were computed for varying degrees of polynomial feedback from the initial state ${\bf x}_0 = [ 10;10;10]$
and simulated to time $T=50$.  The series approximation to the value function and the integral of the running cost
are reported in Table~\ref{table:lorenz}.
\begin{table}[h]
\caption{Lorenz: Value Function Approx.\label{table:lorenz}}
\begin{center}
\begin{tabular}{|c||r@{.}l|r@{.}l|}
\hline
d & \multicolumn{2}{|c|}{$\sum_{i=2}^{d+1} v^{[i]}({\bf x}_0)$} & \multicolumn{2}{|c|}{$\int_0^T \ell({\bf x}(t),{\bf u}(t))dt$} \\
\hline
1 & 7533&49 & 6999&37 \\
2 & 7062&15 & 6911&03 \\
3 & 6957&19 & 6906&45 \\
4 & 6924&27 & 6906&21 \\
5 & 6913&68 & 6906&18 \\
6 & 6910&45 & 6906&17 \\
7 & 6909&30 & 6906&17 \\
\hline
\end{tabular}
\end{center}
\end{table}

\subsection{Ring of van der Pol Oscillators}
As a second test case, we consider controlling a ring of van der Pol oscillators.
\begin{displaymath}
  \ddot{y}_i + (y_i^2-1)\dot{y}_i + y_i = y_{i-1}-2y_i+y_{i+1} + b_i u_i(t),
\end{displaymath} 
for $i=1,\ldots,g$ with $y_i(0)=y_0$ and $\dot{y}_i(0)=0$ (we identify $y_{g+1}=y_1$ and $y_g=y_0$ to 
close the ring).  
The stability of this system
was studied in \cite{nana2006vanderpol} and a related control problem considered in \cite{barron2016vanderpol}.
Choosing different values of $g$ and rewriting as a first-order system of differential equations
allows us to study the cubic-quadratic regulator problem
for problems of size $n=2g$.  We set $b_i$ as 0 or 1 with $m = \| {\bf b} \|_1$.  

In our first experiment, we chose $g=4$, set $b_1=b_2=1$ and computed feedback laws up to septic (7th degree) polynomial
terms.  The problem parameters are thus, $n=8$, $m=2$, $p=3$, and $d=7$.  We also set $y_0=0.3$ and $T=50$ for this study.  The approximations
to the value function are presented in Table~\ref{table:vander1}.  Two points are immediately obvious.  One is that the even degree feedback coefficients are calculated to exactly zero (the state equation only has odd terms).  The second is that the benefit of the control is quickly under the integration threshold by the cubic feedback terms.  The polynomial estimate of the value function slowly continues to improve, but is sufficient by the septic terms.

\begin{table}[h]
\caption{van der Pol: Value Function Approx.\label{table:vander1}}
\begin{center}
\begin{tabular}{|c||r@{.}l|r@{.}l|}
\hline
d & \multicolumn{2}{|c|}{$\sum_{i=2}^{d+1} v^{[i]}({\bf x}_0)$} & \multicolumn{2}{|c|}{$\int_0^T \ell({\bf x}(t),{\bf u}(t))dt$} \\
\hline
1 & 4&6380 & 4&4253 \\
2 & 4&6380 & 4&4253 \\
3 & 4&4125 & 4&4208 \\
4 & 4&4125 & 4&4208 \\
5 & 4&4246 & 4&4208 \\
6 & 4&4246 & 4&4208 \\
7 & 4&4242 & 4&4208 \\
\hline
\end{tabular}
\end{center}
\end{table}

In a second experiment, we used 8 oscillators with 2 controls at nodes 1 and 2.  This led to a control
problem where ${\bf v}_2$ and ${\bf k}_1$ were well defined, ${\bf A}_c$ was stable, yet the origin was only locally stable for the nonlinear system.  Implementing the cubic and quintic controls lead to finite-time blowup.  However, for this same
scenario, choosing $y_0=0.03$ lead to the expected improvements (about 2\%) with higher degree feedback laws
and approximations to the value function that verified the results, see Table~\ref{table:vander2}.
\begin{table}[h]
\caption{van der Pol: Value Function Approx.\label{table:vander2}}
\begin{center}
\begin{tabular}{|c||r@{.}l|r@{.}l|}
\hline
d & \multicolumn{2}{|c|}{$\sum_{i=2}^{d+1} v^{[i]}({\bf x}_0)$} & \multicolumn{2}{|c|}{$\int_0^T \ell({\bf x}(t),{\bf u}(t))dt$} \\
\hline
1 &16&8514 & 16&4579\\
3 &16&0162 & 16&0622\\
5 &16&0830 & 16&0566\\
\hline
\end{tabular}
\end{center}
\end{table}

As a final experiment, we increase the number of controls to 4, but consider different locations for the actuators.
We note that locations at nodes $(1,3,5,7)$ and, by rotational symmetry, $(2,4,6,8)$ lead to uncontrollable $({\bf A},{\bf B})$
pairs.  The ARE must have a solution before any higher degree feedback approximations are defined.  At other locations,
we found success from the original $y_0=0.3$ value.  In Table~\ref{table:vander3}, we list the actuated nodes along with
values for the integrated running cost for linear, cubic, and quintic feedback laws.  In every case except the first, 
there was improvement in the actual performance for higher degree feedback laws.  Some cases showed a 3\% improvement
while others were not significant (0.08\%).  Most of the performance gain was achieved with the addition of the ${\bf k}_3$
term, although there was one example where the cubic feedback lead to finite-time blowup of the solution and 
the full quintic feedback was required to see the performance gains.  The approximation of the value function wasn't as
insightful in deciding when the approximation of the nonlinear feedback gains were of a high enough degree.  A closed-loop
simulation was required to evaluate the performance.
 \begin{table}[h]
\caption{van der Pol: Value Function Approx.\label{table:vander3}}
\begin{center}
\begin{tabular}{|c||r@{.}l|r@{.}l|r@{.}l|}
\hline
nodes & \multicolumn{2}{|c|}{linear} & \multicolumn{2}{|c|}{cubic} & \multicolumn{2}{|c|}{quintic} \\
\hline
(1,2,3,4) & 77&9977 & \multicolumn{2}{|c|}{blow-up} & 75&7120\\
(1,2,3,5) & 29&9355 & 29&1139 & 29&0181\\
(1,2,3,6) &  8&3986 &  8&3910 &  8&3910\\
(1,2,4,5) & 29&4803 & 28&6854 & 28&5952\\
(1,2,4,6) &  7&7364 &  7&7293 &  7&7292\\
(1,2,4,7) &  6&9549 &  6&9489 &  6&9489\\
(1,2,5,6) &  8&8505 &  8&8417 &  8&8417\\
\hline
\end{tabular}
\end{center}
\end{table}

\subsection{Burgers Equation With Reaction Term\label{sec:burgers}}
As a more structured test problem, we consider the PQR problem with a discretization
of the Burgers equation.  This test problem has a long history in the study of control for 
distributed parameter systems, e.g.~\cite{thevenet2009nonlinear}, including the 
development of effective computational methods, e.g.~\cite{burns1990control}.   

We consider the specific problem found in \cite{borggaard2018computation} but with
three control inputs ($m=3$) that consist of uniformly distributed sources over 
disjoint patches.  Thus, we have a bounded input operator.  The formal description of
the problem is 
\begin{displaymath}
  \min_{\bf u}J(z,u) = \int_0^\infty \left(\int_0^1 z^2(\xi,t)\ d\xi + {\bf u}^T(t){\bf u}(t)\right)\ dt
\end{displaymath}
subject to
\begin{align*}
  \dot{z}(x,t) &= \epsilon z_{xx}(x,t) - \frac{1}{2}\left(z^2(x,t)\right)_x  + \alpha z(x,t)
  \\
  &\hspace{.4in}+ \sum_{k=1}^m 
  \chi_{[(k-1)/m,\ k/m]}(x) u_k(t),\\
  z(\cdot,0) &= z_0(\cdot) \in H_{\rm per}^1(0,1),
\end{align*}
where $\chi_{[a,b]}(x)$ is the characteristic function over $[a,b]$. We discretized the state equations with $n$ linear finite elements, set $m=3$, and chose $\epsilon=0.005$ to make the nonlinearity significant and $\alpha=0.3$ to accentuate the need for control in this problem.  

The discretized system fits within the PQR framework (\ref{eq:oc})-(\ref{eq:full}).  The matrices ${\bf A}$, ${\bf B}$ and ${\bf N}_2$ come from the finite element approximation.  The matrix ${\bf Q}_2$ is the finite element mass matrix and the matrix ${\bf R}_2=10~{\bf I}_m$.  For this test, we started with the smooth initial condition
\begin{displaymath}
  z_0(x) = \left\{ \begin{array}{cl}
  0.5 \sin(2\pi x)^2 & x\in(0,0.5) \\
  0 & \mbox{otherwise}
  \end{array} \right. .
\end{displaymath}
Our discretization was performed with 16 linear finite elements.  Approximations to the value function by polynomial expansion and closed-loop numerical simulation to $T=200$ are shown in Table~\ref{table:burgers} (results from a 20 linear element computation are in parenthesis).  We should note that the $n=20$, $m=3$, $p=5$ study reported here, took less that 160 seconds on a 2017 MacBook Pro
with 16GB of RAM.  Therefore, the available Matlab software is sufficiently efficient to study nonlinear feedback on modest sized problems with enough available RAM.

\begin{table}[h]
\caption{Burgers: Value Function Approx.\label{table:burgers}}
\begin{center}
\begin{tabular}{|c||r@{.}l|r@{.}l|}
\hline
d & \multicolumn{2}{|c|}{$\sum_{i=2}^{d+1} v^{[i]}({\bf x}_0)$} & \multicolumn{2}{|c|}{$\int_0^T \ell({\bf x}(t),{\bf u}(t))dt$} \\
\hline
1 & 0&0162721 (0.0175278) & 0&0190134 (0.0190637)\\
2 & 0&0216261 (0.0226721) & 0&0188797 (0.0189653)\\
3 & 0&0200150 (0.0194921) & 0&0187951 (0.0188268)\\
4 & 0&0178709 (0.0172666) & 0&0187623 (0.0188218)\\
5 & 0&0183326 (0.0184925) & 0&0187435 (0.0187726)\\
\hline
\end{tabular}
\end{center}
\end{table}

\section{Conclusions and Future Work}
We presented a special formulation of the Al'Brekht polynomial approximation for
polynomial-quadratic regulator problems.  Writing the system and expansions in terms of Kronecker
products leads to a series of progressively larger linear systems for the next terms
in the expansion.  While easy to write down and implement, efficiency is only achieved
by exploiting new numerical linear algebra tools that avoid the assembly of the large,
dense systems~\citep{kolda2009tensor,chen2019recursive}.  In our previous work~\citep{borggaard2020qqr}, we performed a comparison with a general, well-developed software tool, the Nonlinear Systems
Toolbox~\citep{krener2015NST}, to verify our implementation in the quadratic-quadratic case.  Our solution method was 
competitive with NST in terms of CPU time even if we neglect the overhead in using Matlab's
symbolic toolbox (we described an effective means to compute the derivatives of the
system that are required by NST using automatic differentiation in a previous 
paper~\cite{borggaard2018computation}).

A natural path forward will be to include useful generalizations within our software framework.
This includes the investigation of more general control costs, general mixed state and control terms, addition of descriptor 
systems~\citep{xu1993hamilton}, and the related observer problem.

Finally, we will apply this to more significant applications than the one-dimensional Burgers
equation.  In particular, study how this work could be used in conjunction with reduced
models of complex flows that result in quadratic-in-state systems.

This software is available for download at 
\begin{quote}
https://github.com/jborggaard/QQR
\end{quote}
\section{Acknowledgment}
This research was partially supported by the National Science Foundation under contract DMS-1819110 and the authors gratefully acknowledge the support of the Institute for Mathematics and its Applications (IMA), where this work was initiated during its 
annual program on Control Theory and its Applications.  We are also grateful for the detailed comments of the three reviewers that improved the clarity and completeness of this presentation and informed us of the relevant technical report by~\cite{breiten2017taylor}.

\bibliography{refs}

\end{document}